\documentclass{amsart}
\usepackage{amsfonts}
\usepackage{amssymb}
\usepackage[latin1]{inputenc}
\newtheorem{thm}{Theorem}[section]

\theoremstyle{definition}

\theoremstyle{remark}

\numberwithin{equation}{section}

\numberwithin{equation}{section}

\begin{document}

\title[Conn's linearization theorem for analytic Poisson structures]
{A geometric proof of Conn's linearization theorem for analytic Poisson structures}

\author{Nguyen Tien Zung}
\address{New address starting 01/Sep/2002: Laboratoire Emile Picard,
Mathématiques, Université Toulouse III}
\email{tienzung@math.univ-montp2.fr}

\date{First draft, July 2002}

\begin{abstract}

We give a geometric proof of Conn's linearization theorem for analytic Poisson
structures, without using the fast convergence method.

\end{abstract}
\maketitle

\section{Introduction}

In \cite{Conn-Analytic1984}, Jack Conn proved the following theorem:

\begin{thm}[Conn]
\label{thm:Conn1} Let $\Pi = \Pi_1 + ...$ be an analytic Poisson structure in a
neighborhood of $0$ in ${\mathbb K}^n$ (where $\mathbb K = \mathbb R$ or $\mathbb
C$), which vanishes at $0$ and whose linear part $\Pi_1$ corresponds to a
semi-simple Lie algebra. Then $\Pi$ admits a local analytic linearization at $0$.
\end{thm}

This is probably the first major result about analytic local normal forms for
Poisson structures. A similar result in the smooth case is also obtained by Conn in
\cite{Conn-Smooth1985}.

Conn's proof is based on the powerful method of fast convergence. Using this same
method, in \cite{Zung-Levi2002,MoZu-Levi2002} we extended Conn's results, proving
the Levi decomposition for Poisson structures, and the linearization of Lie
algebroids with a semisimple linear part.

The method of fast convergence has its drawbacks: besides the fact that it often
requires heavy analytical estimations, it also hides the geometrical
picture/structure of the problem. Because of this, Alan Weinstein
\cite{Weinstein-Review1998,Weinstein-Linearization2000} and other people were asking
for a more geometric proof of Theorem \ref{thm:Conn1} and similar results.

The aim of this Note is to sketch such a geometric proof of Theorem \ref{thm:Conn1}.
The techniques used here are first developed in \cite{Zung-Birkhoff2002} where we
studied Birkhoff normal forms for Hamiltonian systems. Our proof consists of the
following steps:

\underline{Step 1}. Find $m$ Casimir functions $F_1,...,F_m$ for the Poisson
structure $\Pi$ in question, where $m$ is the rank of the semisimple Lie algebra
$\frak g$ associated to the linear part $\Pi_1$ of $\Pi$. In order to define these
functions, we will use period integrals (of the symplectic form over 2-cycles in
symplectic leaves) and arguments from \cite{Zung-Birkhoff2002}.

\underline{Step 2}. Use the above Casimir functions to show that the symplectic
foliation (by symplectic leaves) of $\Pi$ is locally analytically diffeomorphic to
the symplectic foliation of $\Pi_1$. So we may assume that $\Pi$ has the same
foliation as $\Pi_1$. In particular, the symplectic foliation of $\Pi$, considered
as a singular foliation (without the symplectic structure), is invariant under the
coadjoint action of $\frak g$.

\underline{Step 3}. Use averaging and Moser's path method to show that there is a
non-autonomous vector field tangent to the symplectic foliation whose time-1 map
moves $\Pi$ to a $\frak g$-invariant Poisson structure. It implies that $\Pi$ admits
a Hamiltonian $\frak g$-action. The components of the corresponding equivariant
moment map will then form a local linear system of coordinates for $\Pi$.

In the next three sections, we will carry out the above three steps, for holomorphic
Poisson structures. Then in Section 5 we will indicate why things work the same in
the real analytic case.

\section{Step 1: Casimir functions}

First let us look at Camisir functions for the linear Poisson structure $\Pi_1$ on
${\frak g}^\ast$, where $\frak g$ is the complex semi-simple Lie algebra of rank $m$
associated to $\Pi_1$. (We will consider $\Pi$ as a Poisson structure in a
neighborhood of $0$ in ${\frak g}^\ast$). It is well known that the symplectic
foliation of $\Pi_1$ (by coadjoint orbits) has codimension $m$, and there are $m$
independent homogeneous polynomial Casimir functions. We will refer to ??? for this
and other informations concerning semi-simple Lie algebras used here.

What is important for us here is the fact that Casmimir functions for $\Pi_1$ can be
defined by period integrals as follows:

Let $P$ be a generic (semisimple regular) coadjoint orbit. Then the dimension of
$H_2(P,\mathbb R)$ is $m$. We can choose a basis $\Gamma_1,\dots,\Gamma_m$ of
$H_2(P,\mathbb R)$, which may be represented by real 2-spheres on $P$ (these
2-spheres can be constructed explicitly by and fixing a Cartan algebra and a root
system).

Define the following period integrals:
$$\rho^i_1 = \int_{\Gamma_i} \omega_1$$
where $\omega^1$ denotes the induced symplectic structure of $\Pi_1$ on $P$ (the
so-called Kostant-Kirillov-Souriau symplectic structure).

So to each generic $P$ we can associate $m$ numbers $\rho_1^1,\dots,\rho_1^m$. By
changing $P$, they become Casimir functions near each generic leaf. But they are not
single-valued (holomorphic) Casimir functions for $\Pi_1$ due to the monodromy
problem: By moving round circle $P$ around a singularity, the 2-cycles on $P$ change
after the move. In other words, the locally flat fiber bundle whose fibers are
$H^2(P,{\mathbb R})$ has a nontrivial monodromy. The corresponding monodromy group
is nothing by the Weyl group (Borel's theorem ?). Thus, if we denote by
$S(\rho_1^1,\dots,\rho_1^m)$ the symmetric algebra of $\rho_1^1,\dots,\rho_1^m$,
then they Weyl group acts on it by monodromy. In fact, the set of polynomial
Casimirs functions for $\Pi_1$ coincides with the set
$$S(\rho_1^1,\dots,\rho_1^m)^W$$
of elements in $$S(\rho_1^1,\dots,\rho_1^m)$$ which are invariant under the Weyl
group action. We may choose a basis $F^1_1 = G^1(\rho_1^1,\dots,\rho_1^m), \dots,
F^m_1 = G^m(\rho_1^1,\dots,\rho_1^m)$ of homogeneous polynomial Casimir functions of
$\Pi_1$.

For examples, when ${\frak g} = sl(m+1,{\mathbb C})$, then $\rho_1^1,\dots,\rho_1^m$
can be chosen to be the first $m$ eigenvalues of $(m+1) \times (m+1)$ matrices,
while $F^1_1,\dots,F^m_1$ are nonlinear symmetric functions of the eigenvalues.

Now look at $\Pi$. An important observation is that, since $\Pi$ is formally
equivalent to $\Pi_1$ (Weintein's theorem \cite{Weinstein-Poisson1983}), most
symplectic leaves of $\Pi$ in a sufficiently small neighborhood of $0$ have $m$
independent 2-cycles inherited from $\Pi_1$. More precisely, for any large natural
number $N$ we have $\Pi = \Pi_1^{(N)} + o(N)$, where $o(N)$ means terms of order
greater than $N$, and $\Pi_1^{(N)}$ is a Poisson structure which is locally
analytically equivalent to $\Pi_1$. Then most symplectic leaves of $\Pi_1^{(N)}$ in
an sufficiently small neighborhood of $0$ are ``nearly tangent'' to symplectic
leaves of $\Pi$. Therefore, due to Reeb's stability, 2-spheres that represent
$2$-cycles on most symplectic leaves of $\Pi_1^{(N)}$ can be projected (in a unique
way homotopically) to $2$-cycles on symplectic leaves of $\Pi$ (This is well defined
outside a horn-shaped neighborhood of the singular set of $\Pi$, see
\cite{Zung-Birkhoff2002} for details). Denote these cycles on symplectic leaves of
$\Pi$ again by $\Gamma_1,\dots,\gamma_m$, and define

$$\rho_i = \int_{\Gamma_i}\omega$$
where $\omega$ is the symplectic form induced from $\Pi$, and

$$F^i = G^i(\rho_1,\dots,\rho_m)$$

Then $F^1,...,F^m$ are single-valued Casimir functions for $\Pi$ outside a
horn-shaped neighborhood of the singular set of $\Pi$. Now make $N$ tend to $\infty$
and use arguments from \cite{Zung-Birkhoff2002} to conclude that $F^1,...,F^m$ are
holomorphic Casimir functions for $\Pi$ in a neighborhood of $0$.

\section{Step 2: Symplectic foliation}

Since $\Pi$ is formally equivalent to $\Pi_1$, the $m$-tuple of Casimir functions
$F^1,...,F^m$ for $\Pi$ is formally equivalent to the $m$-tuple of Casimir functions
$F^1_1,...,F^m_1$ for $\Pi_1$. In other words, if we denote by $y = y(x)$ a formal
diffeomorphism which moves $\Pi$ to $\Pi_1$, then we have

$$F^i (x) = F^i_1 (y(x))$$

Now applying Artin's theorem \cite{Artin-Analytic1968} to the system of analytic
equations
$$F^1_1 (y) - F^1(x) = 0, \dots, F^m_1(y) - F^m(x) = 0$$
and the formal
solution $y = y (x)$, we find a local analytic diffeomorphism $z = z(x)$ (which is
tangent to the formal solution $y = y (x)$ up to any desired order) such that

$$F^i (x) = F^i_1 (z(x))$$

By applying this local diffeomorphism, we may assume that
$$F^i (x) = F^i_1 (x),$$
i.e. the Casimir functions for $\Pi$ are the same as the Casimir functions for
$\Pi_1$. It implies that the symplectic foliation for $\Pi$ is locally the same as
the symplectic foliation for $\Pi_1$, i.e. it is given by coadjoint orbits on
${\frak g}^\ast$ (at least in a dense regular part, but then everywhere in a
neighborhood of $0$ by continuation). We may also assume that $\Pi - \Pi_1 = o(N)$
for some natural number $N$ high enough, i.e.  $\Pi$ is tangent to $\Pi_1$ up to
order $N$.

Though we will not use it in the next section, let us mention the following fact:
because the values of $F^1, \dots, F^m$ on each symplectic leaf determines the
cohomological class of the symplectic form on the leaf (via the period integrals),
the cohomological class of the symplectic form induced by $\Pi$ coincides with the
cohomological class of the Kirillov-Kostant-Souriau symplectic form on each
coadjoint orbit. (One may be tempted to use this fact in order to apply Moser's path
method directly to $\Pi$ and $\Pi_1$).

\section{Step 3: $\frak g$-action and isotopy}

Denote by $G$ the compact group whose Lie algebra is the compact form of our
semi-simple algebra $\frak g$. Then $G$ acts on ${\frak g}^\ast$ by coadjoint
action. For each $g \in G$ and coadjoint orbit $P$ we denote by $\omega^P_g = g^\ast
\omega^P_0$ the image of $\omega^P$ under the action of $g$ on $P$, where $\omega^P$
denotes the induced symplectic form of $\Pi$ on $P$.

Define
$$\omega^P_1 = \int_{g \in G} \omega^P_g d \mu$$
where $\mu$ is the Haar measure on $G$. Since $\Pi$ is tangent to $\Pi_1$ up to
order $N$ and $\Pi_1$ is $G$-invariant, it implies that $\omega^P_1$ is
nondegenerate.

Denote by $\Lambda_1$ the Poisson structure whose symplectic leaves are coadjoint
orbits with symplectic form $\omega^P_1$. It is a holomorphic Poisson structure
(because it is holomorphic at least in the regular part of $\Pi$, and the singular
part is of codimension greater than 1, so we can use Hartogs' extension theorem).

For each coadjoint orbit $P$, put $\omega^P_s = (1-s) \omega^P + s \omega^P_1$. Then
we have an analytic 1-dimensional family of Poisson structure $\Lambda_s$ with
induced symplectic structures $\omega^P_s$, which connects $\Pi = \Lambda_0$ with
$\Lambda_1$. One checks that $\Pi_s$ is well-defined (again by Hartogs' theorem)

Now we will find an analytic flow whose type-$s$ map moves $\Pi$ to $\Lambda_s$ by
Moser's path method. The corresponding time-dependent vector field is given by the
following equation:

$$i_{X_s} \omega^P_s = \alpha$$

where $\alpha$ is an 1-form on each symplectic leaf such that $d\alpha = \omega^P_1
- \omega^P$.

The main difficulty here is how to define $\alpha$. For general singular
foliations/fibrations this would be a highly nontrivial problem. (And if we chose
$\Pi_1$ instead of $\Lambda_1$ in the isotopy, it would be more difficult to define
$\alpha)$. However, our situation here is a little bit special because we have a
compact group action and can define $\alpha$ directly by the following formulas:

For each $g \in G$ denote by $\xi(g) \in \frak g$ the ``smallest'' element of $\frak
g$ such that $g = exp (\xi(g))$ (the set where $\xi(g)$ is not well defined (i.e.
not unique) is of measure 0 so it will not matter), and denote by $X(g)$ the vector
field on ${\frak g}^\ast$ generated by $\xi(g)$. Then we have

$$g^\ast \omega^P - \omega^P = \int_0^1 {\mathcal L}_{X(g)} exp(tX(g))^\ast \omega^P
d t = d (\int_0^1 i_{X(g)} exp(tX(g))^\ast \omega^P d t)$$

and

$$\omega^P_1 - \omega^P = d ( \int_G \int_0^1 i_{X(g)} exp(tX(g))^\ast \omega^P dt d \mu)$$

So we can put

$$\alpha =  \int_G \int_0^1 i_{X(g)} exp(tX(g))^\ast \omega^P dt d \mu$$

This last formula assures the analyticity of the time-dependent vector field $X$
whose time-1 map moves $\Pi$ to $\Lambda_1$.

Applying this analytic time-1 map, we may assume that $\Pi$ is $G$-invariant, which
is the same thing as $\frak g$-invariant. Then since $\frak g$ is semisimple, the
action of $\frak g$ is then Hamiltonian with respect to $\Pi$ and is given by a
(unique) equivariant moment map. Use the components of this moment map as local
coordinates for  a neighborhood of $0$. Then $\Pi$ becomes linear with respect to
these coordinates, and we are done.

\section{The real analytic case}

In the real case, we can still proceed as above. Due to the complex conjugation,
Casimir functions can be chosen to be real. In Step 3, we can still use the compact
group $G$ (which acts in the complex space, not the real one). But due to the
complex conjugation, $\Lambda_1$ and $\alpha$ are real ... \\

Question: does the method presented in this note work in the smooth case? I don't
know, but probably the most difficult part is to show that Casimir functions defined
by period integrals in Step 1 are smooth. For this we probably have to control the
singularity set of $\Pi$ first. Step 2 and Step 3 probably still work, with a few
modifications.


\bibliographystyle{amsplain}

\begin{thebibliography}{1}

\bibitem{Artin-Analytic1968}
M.~Artin, \emph{On the solutions of analytic equations}, Invent. Math.
  \textbf{5} (1968), 277--291.

\bibitem{Conn-Analytic1984}
Jack~F. Conn, \emph{Normal forms for analytic {P}oisson structures}, Ann. of
  Math. (2) \textbf{119} (1984), no.~3, 577--601.

\bibitem{Conn-Smooth1985}
\bysame, \emph{Normal forms for smooth {P}oisson structures}, Ann. of Math. (2)
  \textbf{121} (1985), no.~3, 565--593.

\bibitem{MoZu-Levi2002}
Phillipe Monnier and Nguyen~Tien Zung, \emph{Levi decomposition of smooth
  {P}oisson structures}, in preparation (2002).

\bibitem{Weinstein-Poisson1983}
Alan Weinstein, \emph{The local structure of {P}oisson manifolds}, J.
  Differential Geom. \textbf{18} (1983), no.~3, 523--557.

\bibitem{Weinstein-Review1998}
\bysame, \emph{Poisson geometry}, Differential Geom. Appl. \textbf{9} (1998),
  no.~1-2, 213--238.

\bibitem{Weinstein-Linearization2000}
\bysame, \emph{Linearization problems for {L}ie algebroids and {L}ie
  groupoids}, Lett. Math. Phys. \textbf{52} (2000), no.~1, 93--102.

\bibitem{Zung-Birkhoff2002}
Nguyen~Tien Zung, \emph{Convergence versus integrability in {B}irkhoff normal
  form}, Preprint math.DS/0104279 (2001).

\bibitem{Zung-Levi2002}
\bysame, \emph{Levi decomposition of analytic {P}oisson structures and {L}ie
  algebroids}, Preprint math.DG/0203023 (2002).

\end{thebibliography}

\providecommand{\bysame}{\leavevmode\hbox to3em{\hrulefill}\thinspace}

\end{document}